\documentclass[11pt]{amsart}
\usepackage{latexcad}
\usepackage{natbib}

\def\ldiv{\setminus}
\def\rdiv{/}
\def\text{\textrm}

\theoremstyle{definition}
\newtheorem{example}{Example}       %numbered globally

\theoremstyle{theorem}
\newtheorem{lemma}{Lemma}[section]

\newtheorem{theorem}{Theorem}[section]

\begin{document}

\title[Quasigroups of Bol-Moufang type]{The Varieties of Quasigroups of Bol-Moufang
Type: An Equational Reasoning Approach}

\author{J.~D.~Phillips}

\address{Department of Mathematics \& Computer Science, Wabash College,
Crawfordsville, IN 47933, U.S.A.}

\email{phillipj@wabash.edu}

\author{Petr Vojt\v{e}chovsk\'y}

\address{Department of Mathematics, University of Denver, Denver, CO 80208, U.S.A.}

\email{petr@math.du.edu}

\begin{abstract}
A quasigroup identity is of Bol-Moufang type if two of its three variables
occur once on each side, the third variable occurs twice on each side, the
order in which the variables appear on both sides is the same, and the only
binary operation used is the multiplication, viz. $((xy)x)z=x(y(xz))$. Many
well-known varieties of quasigroups are of Bol-Moufang type. We show that there
are exactly $26$ such varieties, determine all inclusions between them, and
provide all necessary counterexamples. We also determine which of these
varieties consist of loops or one-sided loops, and fully describe the varieties
of commutative quasigroups of Bol-Moufang type. Some of the proofs are
computer-generated.
\end{abstract}



\maketitle

\section{Introduction}

\noindent The purpose of this paper is twofold: to provide the classification
of varieties of quasigroups of Bol-Moufang type, and to demonstrate that the
equational reasoning and finite model builder software currently available is
powerful enough to answer questions of interest in mathematics.

Since we hope to attract the attention of both mathematicians and computer
scientists, we give the necessary background for both groups.

Recall that a class $X$ of universal algebras of the same type is
a \emph{variety} if it is closed under products, subalgebras, and
homomorphic images. Equivalently, $X$ is a variety if and only if
it consists of all universal algebras of the same type satisfying
some identities.

Generally speaking, to establish inclusions between varieties (or sets), it
suffices to use only two types of arguments:
\begin{enumerate}
\item[(i)] Given varieties $\mathcal A$, $\mathcal B$, show that $\mathcal
A\subseteq\mathcal B$.

\item[(ii)] Given varieties $\mathcal A$, $\mathcal B$, decide if there is
$C\in\mathcal B\setminus\mathcal A$.
\end{enumerate}

Throughout this paper, the varieties will be varieties of quasigroups defined
by a single Bol-Moufang identity. Thus, if $i_{\mathcal A}$ is the identity
defining $\mathcal A$, and $i_{\mathcal B}$ is the identity defining $\mathcal
B$, then to establish (i) it suffices to show that $i_{\mathcal A}$ implies
$i_{\mathcal B}$. We use the equational reasoning tool \texttt{Otter} \cite{Mc}
to assist with some of these proofs. As for (ii), we use the finite model
builder \texttt{Mace4} \cite{Mc} to find algebras $C$ in $\mathcal B\setminus
\mathcal A$. We call such algebras $C$ \emph{distinguishing examples}.

The computer-generated \texttt{Otter} proofs are cumbersome and difficult to
read; we do not include them in this paper. However, all have been carefully
organized and are available electronically at \cite{download}. To be able to
read and understand \cite{download}, see Section \ref{Sc:Otter}.

\section{Quasigroups of Bol-Moufang type}

\noindent A set $Q$ with binary operation $\cdot$ is a \emph{quasigroup} if the
equation $a\cdot b=c$ has a unique solution in $Q$ whenever two of the three
elements $a$, $b$, $c\in Q$ are given. Note that multiplication tables of
finite quasigroups are exactly Latin squares.

\begin{figure}[ht]\begin{center}\input{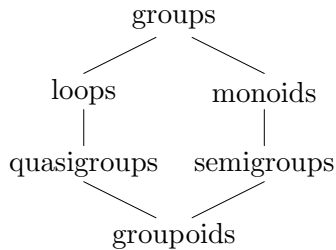}\end{center}
\caption{From groupoids to groups}\label{Fg:Lattice}
\end{figure}

An element $e\in Q$ is called the \emph{left $($right$)$ neutral element} of
$Q$ if $e\cdot a = a$ $(a\cdot e = a)$ holds for every $a\in Q$. An element
$e\in Q$ is the \emph{neutral element} if $e\cdot a = a\cdot e = a$ holds for
every $a\in Q$. In this paper, we use the term \emph{left $($right$)$ loop} for
a quasigroup with a left (right) neutral element. A \emph{loop} is a quasigroup
with a neutral element. Hence loops are precisely `not necessarily associative
groups', as can also be seen from the lattice of varieties depicted in Figure
\ref{Fg:Lattice}. (Recall that a \emph{semigroup} is an associative groupoid,
and a \emph{monoid} is a semigroup with a neutral element.)

%FIGURE
\setlength{\unitlength}{1.2mm}

The above definition of a quasigroup cannot be written in terms of identities,
as it involves existential quantifiers. Fortunately, as is the tradition in the
field, it is possible to introduce a certain kind of universal algebra with 3
binary operations axiomatized by identities (hence forming a variety) that
describes the same objects.

The \emph{variety of quasigroups} consists of universal algebras
$(Q,\,\cdot,\,\ldiv,\,\rdiv)$ whose binary operations $\cdot$, $/$, $\setminus$
satisfy
\begin{displaymath}
    a\cdot(a\ldiv b) = b,\quad (b\rdiv a)\cdot a = b,\quad a\ldiv(a\cdot b) = b,
    \quad (b\cdot a)\rdiv a = b.
\end{displaymath}
It is customary to think of $a\ldiv b$ as division of $b$ by $a$ on the left,
and of $a\rdiv b$ as division of $a$ by $b$ on the right. Note that $a\ldiv b$
is the unique solution to the equation $ax=b$, and similarly for $b\rdiv a$.

This latter description of quasigroups is necessary if one wants to work with
equational reasoning software, such as \texttt{Otter}.

We are coming to the title definition of this paper. An identity $\alpha=\beta$
is of \emph{Bol-Moufang type} if (i) the only operation in $\alpha$, $\beta$ is
$\cdot$, (ii) the same 3 variables appear on both sides, in the same order,
(iii) one of the variables appears twice on both sides, (iv) the remaining two
variables appear once on both sides. For instance, $(x\cdot y)\cdot(x\cdot z) =
x\cdot((y\cdot x)\cdot z)$ is an identity of Bol-Moufang type. A systematic
notation for identities of Bol Moufang type was introduced in \cite{PV}, and
will be reviewed in Section \ref{Sc:SN}.

A variety of quasigroups (loops) is said to be of \emph{Bol-Moufang} type if it
is defined by one additional identity of Bol-Moufang type.

We say that two identities (of Bol-Moufang type) are \emph{equivalent} if they
define the same variety. This definition must be understood relative to some
underlying variety, since two identities can be equivalent for loops but not
for quasigroups, as we shall see.

Several well-known varieties of loops are of Bol-Moufang type. Their
classification was initiated by Fenyves in \cite{Fe1}, \cite{Fe2}, and
completed by the authors in \cite{PV}. None of these three papers required
computer calculations. However, shortly after this current project was
undertaken, it became obvious that the situation for varieties of quasigroups
of Bol-Moufang type is more intricate and complex. That is why we opted for
presenting the results in this format, with the lengthier computer-generated
proofs omitted.

\section{Otter and Mace4}\label{Sc:Otter}

\noindent Our investigations were aided by the automated reasoning tool
\texttt{Otter} and by the finite model builder \texttt{Mace4} \cite{Mc}, both
developed by William McCune (Argonne). \texttt{Otter} implements the
Knuth-Bendix algorithm, and has proven to be effective at equational reasoning
\cite{Mc1}. See \cite{MP} for more about \texttt{Otter}'s technical
specifications, as well as links to results assisted by \texttt{Otter} and
\texttt{Mace4}. A self-contained, elementary introduction to \texttt{Otter} can
be found in \cite{Ph}.

Many authors simply use the \texttt{Otter} output file as the proof of a
theorem; it is common practice to publish untranslated \texttt{Otter} proofs
\cite{MP}. This is mathematically sound since the program can be made to output
a simple \emph{proof object}, which can be independently verified by a short
\texttt{lisp} program. We have posted all proofs omitted in this paper at
\cite{download}. In the proofs at \cite{download}, ``para\_into'' and
``para\_from'' are short for ``paramodulation into'' and ``paramodulation
from'', and they are the key steps in any \texttt{Otter} proof. Very crudely,
paramodulation is an inference rule that combines variable instantiation and
equality substitution into one step \cite{MP}.

The proofs generated by \texttt{Otter} contain all information necessary for
their translation into human language; nevertheless, they are not easy to read.
This is because \texttt{Otter} often performs several nontrivial substitutions
at once. Many of the proofs can be made significantly shorter, especially with
some knowledge of the subject available, however, some proofs appear to be
rather clever even after being translated. In other words, one often obtains no
insight into the problem upon seing the \texttt{Otter} proof. We would be happy
to see more intuitive proofs, but we did not feel that they are necessary for
our programme. If the reader wants to come up with such proofs, he/she should
be aware of the standard techniques of the field, such as autotopisms,
pseudo-automorphisms, and their calculus \cite{Br}, \cite{Pf}.

\texttt{Mace4} is a typical finite model builder. Thus, given a finite set of
identities (or their negations), it attempts to construct a universal algebra
satisfying all of the identities. Given the huge number of nonisomorphic (or
nonisotopic) quasigroups of even small orders (cf. \cite[p.\ 61]{Pf}), it is
not easy to construct such examples by hand, without some theory.
\texttt{Mace4} was therefore invaluable for the purposes of this work,
specifically Section \ref{Sc:DE}.

\section{Systematic notation}\label{Sc:SN}

\noindent The following notational conventions will be used throughout the
paper. We omit $\cdot$ while multiplying two elements (eg $x\cdot y = xy$), and
reserve $\cdot$ to indicate priority of multiplication (eg $x\cdot yz =
x(yz)$). Also, we declare $\ldiv$ and $\rdiv$ to be less binding than the
omitted multiplication (eg $x\rdiv yz = x\rdiv (yz)$), and if $\cdot$ is used,
we consider it to be less binding than any other operation (eg $x\cdot yz\ldiv
y = x((yz)\ldiv y)$).

Let $x$, $y$, $z$ be all the variables appearing in an identity of Bol-Moufang
type. Without loss of generality, we can assume that they appear in the terms
in alphabetical order. Then there are exactly $6$ ways in which the $3$
variables can form a word of length $4$, and there are exactly $5$ ways in
which a word of length $4$ can be bracketed, namely:
\begin{displaymath}
\begin{array}{cc}
    \begin{array}{c|c}
        A&xxyz\\
        B&xyxz\\
        C&xyyz\\
        D&xyzx\\
        E&xyzy\\
        F&xyzz
    \end{array}
    \quad\quad&
    \begin{array}{c|c}
        1&o(o(oo))\\
        2&o((oo)o)\\
        3&(oo)(oo)\\
        4&(o(oo))o\\
        5&((oo)o)o
    \end{array}
\end{array}
\end{displaymath}
Let $Xij$ with $X\in\{A,\dots,F\}$, $1\le i<j\le 5$ be the identity whose
variables are ordered according to $X$, whose left-hand side is bracketed
according to $i$, and whose right-hand side is bracketed according to $j$. For
instance, $C25$ is the identity $x((yy)z)=((xy)y)z$.

It is now clear that any identity of Bol-Moufang type can be transformed into
some identity $Xij$ by renaming the variables and interchanging the left-hand
side with the right-hand side. There are therefore $6\cdot(4+3+2+1)=60$
``different'' identities of Bol-Moufang type, as noted already in \cite{Fe2},
\cite{Ku2}, \cite{PV}.

The \emph{dual} of an identity $I$ is the identity obtained from $I$ by reading
it backwards, i.e., from right to left. For instance, the dual of
$(xy)(xz)=((xy)x)z$ is the identity $z(x(yx)) = (zx)(yx)$. With the above
conventions in mind, we can rewrite the latter identity as $x(y(zy)) =
(xy)(zy)$. One can therefore identify the dual of any identity $Xij$ with some
identity $(Xij)'=X'j'i'$. The name $X'j'i'$ of the dual of $Xij$ is easily
calculated with the help of the following rules:
\begin{displaymath}
    A'=F,\quad B'=E,\quad C'=C,\quad D'=D,\quad 1'=5,\quad 2'=4,\quad 3'=3.
\end{displaymath}

\section{Canonical definitions of some varieties of quasigroups}

\noindent Table \ref{Tb:Definitions} defines $26$ varieties of quasigroups. As
we shall see, these varieties form the complete list of quasigroup varieties of
Bol-Moufang type.

\begin{table}\caption{Definitions of varieties of quasigroups}\label{Tb:Definitions}
\begin{small}
\begin{displaymath}
    \begin{array}{ccccc}
        \text{variety}&\text{abbrev.}&\text{defining\
        identity}&\text{its\ name}&\text{reference}\\
        \hline
        \text{groups}&
            \text{GR}&x(yz)=(xy)z&&\text{folklore}\\
       \text{RG1-quasigroups}&
            \text{RG1}&x((xy)z)=((xx)y)z&A25&\text{new}\\
        \text{LG1-quasigroups}&
            \text{LG1}&x(y(zz))=(x(yz))z&F14&\text{new}\\
        \text{RG2-quasigroups}&
            \text{RG2}&x(x(yz))=(xx)(yz)&A23&\text{new}\\
        \text{LG2-quasigroups}&
            \text{LG2}&(xy)(zz)=(x(yz))z&F34&\text{new}\\
        \text{RG3-quasigroups}&
            \text{RG3}&x((yx)z)=((xy)x)z&B25&\text{new}\\
        \text{LG3-quasigroups}&
            \text{LG3}&x(y(zy))=(x(yz))y&E14&\text{new}\\
        \text{extra q.}&
           \text{EQ}&x(y(zx))=((xy)z)x&D15&\text{\cite{Fe1}}\\
        \text{Moufang q.}&
           \text{MQ}&(xy)(zx)=(x(yz))x&D34&\text{\cite{Mo}}\\
        \text{left Bol q.}&
           \text{LBQ}&x(y(xz))=(x(yx))z&B14&\text{\cite{Ro}}\\
        \text{right Bol q.}&
           \text{RBQ}&x((yz)y)=((xy)z)y&E25&\text{\cite{Ro}}\\
        \text{C-quasigroups}&
            \text{CQ}&x(y(yz))=((xy)y)z&C15&\text{\cite{Fe2}}\\
        \text{LC1-quasigroups}&
            \text{LC1}&(xx)(yz)=(x(xy))z&A34&\text{\cite{Fe2}}\\
        \text{LC2-quasigroups}&
            \text{LC2}&x(x(yz))=(x(xy))z&A14&\text{new}\\
        \text{LC3-quasigroups}&
            \text{LC3}&x(x(yz))=((xx)y)z&A15&\text{new}\\
        \text{LC4-quasigroups}&
            \text{LC4}&x(y(yz))=(x(yy))z&C14&\text{new}\\
        \text{RC1-quasigroups}&
            \text{RC1}&x((yz)z)=(xy)(zz)&F23&\text{\cite{Fe2}}\\
        \text{RC2-quasigroups}&
            \text{RC2}&x((yz)z)=((xy)z)z&F25&\text{new}\\
        \text{RC3-quasigroups}&
            \text{RC3}&x(y(zz))=((xy)z)z&F15&\text{new}\\
        \text{RC4-quasigroups}&
            \text{RC4}&x((yy)z)=((xy)y)z&C25&\text{new}\\
        \text{left alternative q.}&
            \text{LAQ}&x(xy)=(xx)y&&\text{folklore}\\
        \text{right alternative q.}&
            \text{RAQ}&x(yy)=(xy)y&&\text{folklore}\\
        \text{flexible q.}&
            \text{FQ}&x(yx)=(xy)x&&\text{\cite{Pf}}\\
        \text{left nuclear square q.}&
            \text{LNQ}&(xx)(yz)=((xx)y)z&A35&\text{new}\\
        \text{middle nuclear square q.}&
            \text{MNQ}&x((yy)z)=(x(yy))z&C24&\text{new}\\
        \text{right nuclear square q.}&
            \text{RNQ}&x(y(zz))=(xy)(zz)&F13&\text{new}\\
    \end{array}
\end{displaymath}
\end{small}
\end{table}

We have carefully chosen the defining identities in such a way that they are
either self-dual (GR, EQ, CQ, FQ, MNQ) or coupled into dual pairs (Lx${}'$=Rx).
The only exception to this rule is the Moufang identity $D34$. We will often
appeal to this duality.

The reasoning behind the names of the new varieties in Table
\ref{Tb:Definitions} is as follows: If any of the quasigroups LGi, RGi is a
loop, it becomes a group; if any of the quasigroups LCi, RCi is a loop, it
becomes an LC-loop, RC-loop, respectively. All of this will be clarified in the
next section.

Although we will use the new names and abbreviations of Table
\ref{Tb:Definitions} for the rest of the paper, the reader is warned that other
names exist in the literature. Fenyves \cite{Fe2} assigned numbers to the $60$
identities of Bol-Moufang type in a somewhat random way, and Kunen \cite{Ku2}
developed a different systematic notation that does not reflect inclusions
between varieties. Our notation is an extension of \cite{PV}.

\section{Loops of Bol-Moufang type}

\begin{figure}[ht]\begin{center}\input{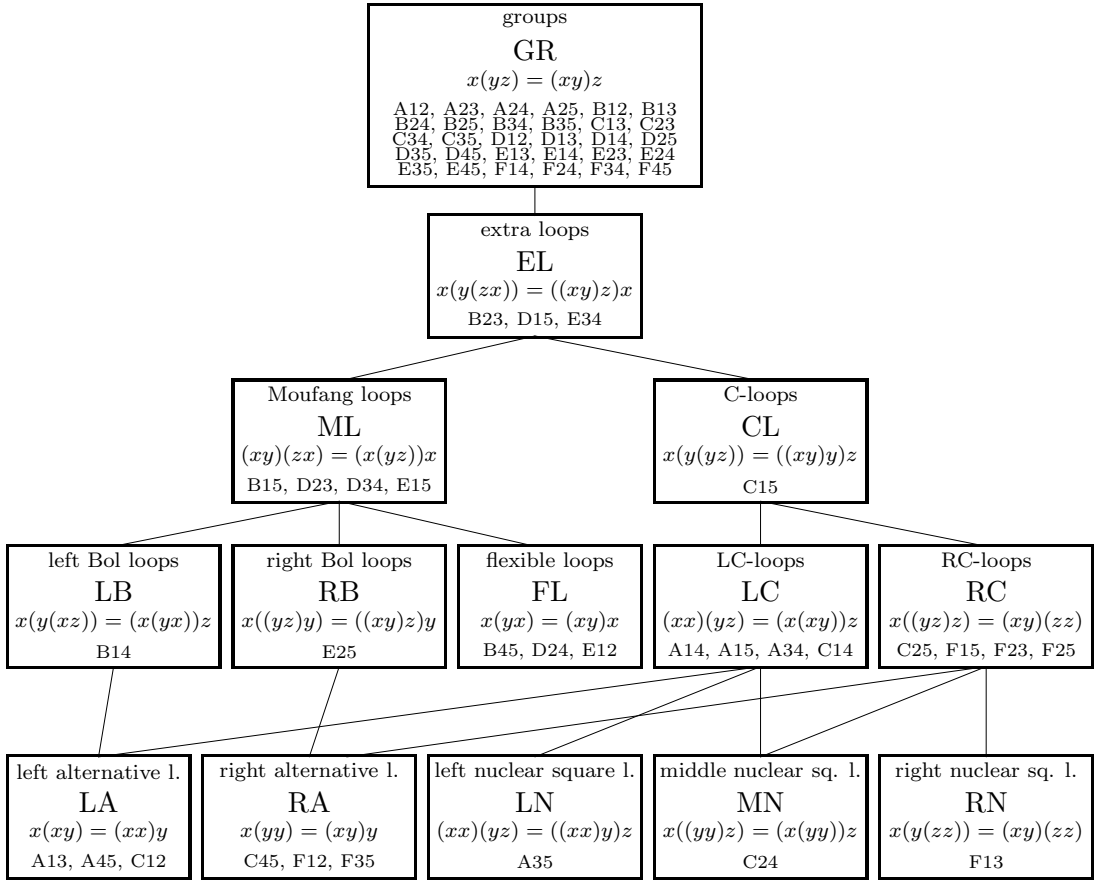}\end{center}
\caption{Varieties of loops of Bol-Moufang type}\label{Fg:LTree}
\end{figure}

\noindent The varieties of loops of Bol-Moufang type were fully described in
\cite{PV}. The situation is summarized in Figure \ref{Fg:LTree}---a Hasse
diagram of varieties of loops of Bol-Moufang type. For every variety in the
diagram we give: the name of the variety, abbreviation of the name, defining
identity, and all equivalent Bol-Moufang identities defining the variety.
Inclusions among varieties are indicated by their relative position and
connecting lines, as is usual in a Hasse diagram. The higher a variety is in
Figure \ref{Fg:LTree}, the smaller it is.

We use Figure \ref{Fg:LTree} as the starting point for our investigation of the
quasigroup case.

%FIGURE
\setlength{\unitlength}{1mm}

\section{Equivalences}

\noindent Our first task is to determine which quasigroup varieties of
Bol-Moufang type are equivalent. Clearly, if two identities are equivalent for
quasigroups, they must also be equivalent for loops, but not \emph{vice versa}.

The results can be found in Figure \ref{Fg:QTree}. (Note that the right half of
Figure \ref{Fg:QTree} is dual to the left side, with the exception of CQ. The
$5$ varieties in the middle and CQ are self-dual.) Actually, all sets of
equivalent identities in Figure \ref{Fg:QTree} are maximal equivalent sets, but
we will not verify the maximality yet. (See Section \ref{Sc:Ids}.) Let us go
through diagram \ref{Fg:QTree} and comment on it.

%FIGURE
\setlength{\unitlength}{0.96mm}
\begin{figure}[ht]\centering\input{qdiag2.lp}
\caption{Varieties of quasigroups of Bol-Moufang type}\label{Fg:QTree}
\end{figure}

We note that there is nothing to show if we list a single identity for a given
variety in Figure \ref{Fg:QTree}, since we merely need to show that all
identities in a given box of Figure \ref{Fg:QTree} are equivalent. This takes
care of $18$ varieties.

Kunen \cite{Ku} showed that the identities $B15$, $D23$, $D34$, $E15$ are
equivalent for quasigroups, hence settling MQ. He also showed that the three EQ
identities $B23$, $D15$, $E34$ are equivalent \cite{Ku2}. We will talk more
about these results in Section \ref{Sc:Ids}.

\begin{lemma} All identities listed under \emph{GR} in Figure $\ref{Fg:QTree}$ are
equivalent, and imply associativity.
\end{lemma}
\begin{proof}
It is very easy to show, just as we did in \cite[Proposition 4.1]{PV}, that all
the identities with the exception of $C23$ and $(C23)'=C34$ imply
associativity. (To show that $C23$ implies associativity, we took advantage of
a neutral element in \cite{PV}.)

Assume that a quasigroup $Q$ satisfies $C23$. We will first show that $Q$
possesses a right neutral element.

Substitute $(yy)\ldiv z$ for $z$ in $C23$ to obtain
\begin{equation}\label{Eq:Aux1}
    xz = x(yy\cdot yy\ldiv z) = (xy)(y\cdot yy\ldiv z).
\end{equation}
With $x\rdiv y$ instead of $x$ and $y$ instead of $z$ we obtain
\begin{displaymath}
    x = (x\rdiv y)y = (x\rdiv y\cdot y)(y\cdot yy\ldiv y) = x(y\cdot yy\ldiv y).
\end{displaymath}
This means that $y(yy\ldiv y)=e$ is independent of $y$, and is the right
neutral element of $Q$. Since $x(x\ldiv x)=x$, too, the right neutral element
can be written as $e=x\ldiv x = y\ldiv y$.

Now, $(\ref{Eq:Aux1})$ implies
\begin{equation}\label{Eq:Aux2}
    xy\ldiv xz = y\cdot yy\ldiv z,
\end{equation}
which, with $x\ldiv y$ instead of $y$, implies
\begin{equation}\label{Eq:Aux3}
    y\ldiv xz = x\ldiv y\cdot (x\ldiv y)(x\ldiv y)\ldiv z.
\end{equation}
Since $xe=x$, we can use $e$ instead of $z$ in $(\ref{Eq:Aux2})$ to get
\begin{equation}\label{Eq:Aux4}
    xy\ldiv x = y(yy\ldiv e).
\end{equation}
Finally, from $e=y(yy\ldiv y)$, we have
\begin{equation}\label{Eq:Aux5}
    yy\ldiv y = y\ldiv e.
\end{equation}
We now show that $e\ldiv y = y$, which is equivalent to saying that $y=ey$, or,
that $e$ is the neutral element of $Q$. The associativity of the loop $Q$ then
easily follows, just as in \cite{PV}.

We have $e\ldiv y = y(y\ldiv e)\ldiv y$, by one of the quasigroup axioms. Using
$(\ref{Eq:Aux4})$ with $y$ instead of $x$ and $y\ldiv e$ instead of $y$, we can
write
\begin{displaymath}
    e\ldiv y = y(y\ldiv e)\ldiv y =
    y\ldiv e\cdot (y\ldiv e)(y\ldiv e)\ldiv e.
\end{displaymath}
By $(\ref{Eq:Aux5})$, the right-hand side can be written as
\begin{displaymath}
    yy\ldiv y\cdot (yy\ldiv y)(yy\ldiv y)\ldiv e.
\end{displaymath}
By $(\ref{Eq:Aux3})$ with $yy$ instead of $x$ and $e$ instead of $z$, this is
equal to $y\ldiv(yy)e=y\ldiv yy=y$, and we are through.
\end{proof}

\begin{lemma} Each of the three \emph{FQ} quasigroup identities $B45$, $D24$, $E12$ is
equivalent to flexibility.
\end{lemma}
\begin{proof} By duality, it suffices to show that $B45$ is equivalent
to flexibility, and that $D24$ is equivalent to flexibility. Cancel $z$ on the
right in $B45$ to obtain flexibility. Substituting $u$ for $yz$ in $D24$ yields
flexibility. Flexibility clearly implies both $B45$ and $D24$.
\end{proof}

\begin{lemma}
The \emph{LG1} identities $D14$ and $F14$ are equivalent for quasigroups.
\end{lemma}
\begin{proof}
By \texttt{Otter}, \cite[Thms. 1, 2]{download}.
\end{proof}

\begin{lemma}
Each of the \emph{LAQ} identities $A13$, $A45$ and $C12$ is equivalent to the
left alternative law.
\end{lemma}
\begin{proof}
It is easy to see that any of these three identities implies the left
alternative law (e.g. let $u=yz$ in $A13$, cancel $z$ on the right in $A45$,
cancel $x$ on the left in $C12$). With the left alternative law available, all
three identities clearly hold.
\end{proof}

\section{Neutral elements in quasigroups of Bol-Moufang type}
\label{Sc:Ids}

\noindent In 1996, Kunen \cite{Ku} discovered the surprising fact that a
Moufang quasigroup is necessarily a (Moufang) loop. In other words, if a
quasigroup $Q$ satisfies any of the Moufang identities (cf. MQ in Figure
\ref{Fg:QTree}), it possesses a neutral element. He then went on to investigate
this property for many quasigroups of Bol-Moufang type \cite{Ku2}.

Since we wish to translate the diagram \ref{Fg:LTree} into the quasigroup case,
we need to know which of the quasigroup varieties defined in Table
\ref{Tb:Definitions} are in fact loop varieties. It will also be useful later
to know if the varieties consist of one-sided loops.

The results are summarized in Figure \ref{Fg:QTree}. The superscript
immediately following the abbreviation of the name of the variety $\mathcal A$
indicates if:
\begin{enumerate}
\item[(2)] every quasigroup in $\mathcal A$ is a loop,

\item[(L)] every quasigroup in $\mathcal A$ is a left loop, and there is
$Q\in\mathcal A$ that is not a loop,

\item[(R)] every quasigroup in $\mathcal A$ is a right loop, and there is
$Q\in\mathcal A$ that is not a loop,

\item[(0)] there is $Q_L\in\mathcal A$ that is not a left loop, and there is
$Q_R\in\mathcal A$ that is not a right loop.
\end{enumerate}
The justification follows.

\subsection{Neutral elements}

GR consists of groups, hence loops. MQ consists of loops by the above-mentioned
result of Kunen \cite{Ku}. EQ consists of loops by \cite{Ku2}. Every
MNQ-quasigroup and every LC1-quasigroup is a loop by \cite[Thm. 3.1]{Ku2}.

\subsection{One-sided neutral elements}

If one wants to show that a quasigroup $Q$ is a right loop, it suffices to
prove that $x\ldiv x = y\ldiv y$ for every $x$, $y\in Q$, since then we can set
$e=x\ldiv x$ for some fixed $x$, and we get $y\cdot e = y\cdot x\ldiv x =
y\cdot y\ldiv y = y$ for every $y$. Similarly, $Q$ will be a left loop if
$x\rdiv x = y\rdiv y$ for every $x$, $y\in Q$.

\begin{lemma} Every \emph{LG3}-quasigroup is a right loop. Every
\emph{LNQ}-quasigroup is a left loop. Every \emph{LC4}-quasigroup is a right
loop. Every \emph{LBQ}-quasigroup is a right loop. Every \emph{LAQ}-quasigroup
is a left loop.
\end{lemma}
\begin{proof}
Let $Q$ be an LG3-quasigroup. Fix $x$, $y\in Q$ and choose $z\in Q$ so that
$y(zy)=y$. Then $xy = x(y(zy))=x(yz)\cdot y$ by $E14$, and thus $x = x(yz)$, or
$x\ldiv x = yz$. Since $z$ depends only on $y$, we see that $x\ldiv x$ is
independent of $x$. In other words, $x\ldiv x = y\ldiv y$ for every $x$, $y\in
Q$.

Let $Q$ be an LNQ-quasigroup. Fix $x\in Q$ and choose $y\in Q$ such that
$xx\cdot y = xx$. Then $(xx)(yz) = ((xx)y)z = xx\cdot z$ by $A35$, and thus
$yz=z$, or $z\rdiv z = y$. Since $y$ does not depend on $z$, we are done.

Let $Q$ be an LC4-quasigroup. We have $x = (x\rdiv yy)(yy) = (x\rdiv yy)(y\cdot
y(y\ldiv y))$, which is by $C14$ equal to $(x\rdiv yy\cdot yy)(y\ldiv y) =
x(y\ldiv y)$. Thus $x\ldiv x=y\ldiv y$ and $Q$ is a right loop.

Let $Q$ be an LBQ-quasigroup. Then $xy = x(y\rdiv x \cdot x(x\ldiv x))$, which
is by $B14$ equal to $(x(y\rdiv x \cdot x))(x\ldiv x) = (xy)(x\ldiv x)$. With
$u=xy$, we see that $u=u(x\ldiv x)$, or $u\ldiv u = x\ldiv x$.

Let $Q$ be a quasigroup satisfying the left alternative law $x(xy)=(xx)y$. With
$xx\ldiv y$ instead of $y$ we obtain $x(x\cdot xx\ldiv y) = y$, or
\begin{equation}\label{Eq:LAQ1}
    xx\ldiv y = x\ldiv (x\ldiv y).
\end{equation}
The left alternative law also yields $x(xy)\rdiv y = xx$, which, upon
substituting $x\rdiv x$ for $x$ and $x$ for $y$, becomes $x\rdiv x = x\rdiv
x\cdot x\rdiv x$. Then $(\ref{Eq:LAQ1})$ can be used to conclude that $(x\rdiv
x)\ldiv y = (x\rdiv x \cdot x\rdiv x)\ldiv y = (x\rdiv x)\ldiv ((x\rdiv x)\ldiv
y)$. Upon multiplying this equation on the left by $x\rdiv x$, we get $y =
(x\rdiv x)\ldiv y$, or $x\rdiv x = y\rdiv y$.
\end{proof}

\begin{lemma} Every \emph{LG2}-quasigroup is a right loop. Every
\emph{LC3}-quasigroup is a left loop.
\end{lemma}
\begin{proof}
See \cite[Thm.\ 3]{download} and \cite[Thm.\ 4]{download}.
\end{proof}

\subsection{Missing one-sided neutral elements}

All examples of order $n$ below have elements $0$, $1$, $\dots$, $n-1$, and
their multiplication tables have columns and rows labelled by $0$, $\dots$,
$n-1$, in this order. Most of the Examples \ref{Ex:A1}--\ref{Ex:A5} will also
be used in Section \ref{Sc:Imps}.

\begin{example}\label{Ex:A1}
LNQ-quasigroup that is not a right loop:
\begin{displaymath}
\begin{array}{cccc}
    0&1&2&3\\
    2&0&3&1\\
    1&3&0&2\\
    3&2&1&0
\end{array}
\end{displaymath}
For another example, see $I(2,1,3)$ of \cite{Ku2}.
\end{example}

\begin{example}\label{Ex:A2}
FQ- and LC2-quasigroup that is neither a left loop nor a right loop:
\begin{displaymath}
\begin{array}{ccc}
    0&2&1\\
    2&1&0\\
    1&0&2
\end{array}
\end{displaymath}
\end{example}

\begin{example}\label{Ex:A3}
LG1-, LG2-, LG3- and LBQ-quasigroup that is not a left loop:
\begin{displaymath}
\begin{array}{ccc}
    0&2&1\\
    1&0&2\\
    2&1&0
\end{array}
\end{displaymath}
\end{example}

\begin{example}\label{Ex:A4}
LC3-quasigroup that is not a right loop:
\begin{displaymath}
\begin{array}{ccccc}
    0&1&2&3&4\\
    2&3&1&4&0\\
    3&0&4&2&1\\
    4&2&0&1&3\\
    1&4&3&0&2
\end{array}
\end{displaymath}
\end{example}

\begin{example}\label{Ex:A5}
A left alternative quasigroup that is not a a right loop.
\begin{displaymath}
\begin{array}{cccccc}
    0& 1& 2& 3& 4& 5\\
    2& 3& 4& 0& 5& 1\\
    1& 4& 3& 5& 0& 2\\
    4& 0& 5& 2& 1& 3\\
    3& 5& 0& 1& 2& 4\\
    5& 2& 1& 4& 3& 0
\end{array}
\end{displaymath}
\end{example}

\subsection{Summary for neutral elements}

We have recalled that LC1, GR, EQ, MQ, MNQ and RC1 consist of loops. Moreover,
if the inclusions of Figure \ref{Fg:QTree} are correct (and they are, as we
will see in Section \ref{Sc:Imps}), we have just established that all varieties
in the left half of diagram \ref{Fg:QTree} behave as indicated with respect to
neutral element. The right half of the diagram then follows by duality.

\section{Implications}\label{Sc:Imps}

\noindent Many of the implications in diagram \ref{Fg:QTree} can be established
as follows: Let $\mathcal A$ be a quasigroup variety consisting of loops, let
$\mathcal B\supseteq \mathcal A$ be a loop variety, and let $\mathcal C$ be a
quasigroup variety defined by any of the (equivalent defining) identities of
$\mathcal B$. Then $\mathcal C\supseteq \mathcal B$, and thus $\mathcal
C\supseteq \mathcal A$. For example: LC1 is contained in LNQ because LC1 is a
quasigroup variety consisting of loops by Section \ref{Sc:Ids}, and because
LC1$=$LC$\subseteq$LN as loops, by \cite{PV}.

The four inclusions GR$\subseteq$LG1, GR$\subseteq$LG2, GR$\subseteq$RG2,
GR$\subseteq$RG1 are trivial, as GR implies everything (associativity).

\begin{lemma} Every \emph{LC1}-quasigroup is an \emph{LC3}-quasigroup.
\end{lemma}
\begin{proof} We know that an LC1-quasigroup is left alternative. Then
$x(x\cdot yz)=xx\cdot yz=(x\cdot xy)z=(xx\cdot y)z$, where the middle equality
is just $A34$.
\end{proof}

The remaining 10 inclusions are labelled by asterisk (*) in Figure
\ref{Fg:QTree}, and can be verified by \texttt{Otter} as follows:

An LG1-quasigroup is an LG3-quasigroup by \cite[Thm.\ 5]{download}, an
LG1-quasigroup is an LC4-quasigroup by \cite[Thm.\ 6]{download}, an
LC4-quasigroup is an LC2-quasigroup by \cite[Thm.\ 7]{download}, an
LG1-quasigroup is an LBQ-quasigroup by \cite[Thm.\ 8]{download}, and an
LC1-quasigroup is an LC4-quasigroup by \cite[Thm.\ 9]{download}.

The remaining 5 inclusions follow by duality.

At this point, we have justified all inclusions in Figure \ref{Fg:QTree}. We
have not yet shown that the inclusions are proper (i.e., that the sets of
equivalent identities are maximal), and that no inclusions are missing. All of
this will be done next.

\section{Distinguishing examples}\label{Sc:DE}

\noindent Since there are 26 varieties to be distinguished, we will proceed
systematically. Our strategy is as follows:

There are 17 minimal elements (maximal with respect to inclusion) in Figure
\ref{Fg:QTree}, namely LC2, LG3, LBQ, LC3, LNQ, LG2, LAQ, FQ, MNQ, CQ, RAQ,
RG2, RNQ, RBQ, RC3, RG3 and RC2. Assume we want to find $C\in\mathcal
A\setminus \mathcal B$, for some varieties $\mathcal A$, $\mathcal B$ of Figure
\ref{Fg:QTree}. It then suffices to find $C\in\mathcal A\setminus\mathcal C$,
where $\mathcal C$ is any of the minimal elements below $\mathcal B$ (but, of
course, not below $\mathcal A$). For example, in order to distinguish $LG1$
from the varieties not below LG1 in Figure \ref{Fg:QTree}, we only need
$14=17-3$ distinguishing examples.

Assume that given $\mathcal A$, all distinguishing examples $C\in\mathcal
A\setminus \mathcal B$ were found. Then our task is much simpler for any
$\mathcal C$ below $\mathcal A$, since we only must look at minimal elements
below $\mathcal A$ that are not below $\mathcal C$. For instance, when all
distinguishing examples for LG1 are found, we only need 2 additional
distinguishing examples for LG3, namely $C\in LG3\setminus LC2$ and $D\in
LG3\setminus LBQ$.

We get many distinguishing examples for free by taking advantage of the
existence of (one-sided) neutral elements in the varieties under consideration.
For instance, there must be a quasigroup in LG1$\setminus$LC3, as it suffices
to take any LG1-quasigroup that is not a left loop. Such examples will be
called of \emph{type $1$}.

Some of the remaining distinguishing examples can be obtained from the results
of \cite{PV}. For instance, there must be a quasigroup in LC4$\setminus$LBQ,
because, by \cite{PV}, there is a loop in LC$\setminus$LB. Such examples will
be called of \emph{type $2$}.

Finally, about half of the examples follows by duality. Hence, it suffices to
find all distinguishing examples of the form $C\in\mathcal A\setminus\mathcal
B$, where $\mathcal A$ is in the left half of Figure \ref{Fg:QTree}, including
EQ, MQ, FQ, MNQ and CQ (these five varieties are self-dual, as we have already
noticed).

It turns out that only 9 additional examples are needed for a complete
discussion. Here they are:

\begin{example}\label{Ex:E1}
An LG1-quasigroup with additional properties (see below).
\begin{displaymath}
\begin{array}{cccccc}
    0& 2& 1& 3& 5& 4\\
    1& 4& 0& 5& 3& 2\\
    2& 0& 5& 4& 1& 3\\
    3& 5& 4& 0& 2& 1\\
    4& 1& 3& 2& 0& 5\\
    5& 3& 2& 1& 4& 0
\end{array}
\end{displaymath}
\end{example}

\begin{example}\label{Ex:E5}
An LG1-quasigroup with additional properties (see below).
\begin{displaymath}
\begin{array}{cccccc}
    0& 3& 4& 1& 2& 5\\
    1& 2& 5& 0& 3& 4\\
    2& 1& 0& 5& 4& 3\\
    3& 0& 1& 4& 5& 2\\
    4& 5& 2& 3& 0& 1\\
    5& 4& 3& 2& 1& 0
\end{array}
\end{displaymath}
\end{example}

\begin{example}\label{Ex:E8}
An LC4-quasigroup with additional properties (see below).
\begin{displaymath}
\begin{array}{cccccc}
    0& 1& 3& 2& 5& 4\\
    1& 5& 0& 4& 2& 3\\
    2& 0& 4& 5& 3& 1\\
    3& 4& 5& 0& 1& 2\\
    4& 2& 1& 3& 0& 5\\
    5& 3& 2& 1& 4& 0
\end{array}
\end{displaymath}
\end{example}

\begin{example}\label{Ex:E10}
An LG3-quasigroup with additional properties (see below).
\begin{displaymath}
\begin{array}{cccccc}
    0& 2& 3& 1& 4& 5\\
    1& 0& 4& 5& 2& 3\\
    2& 4& 5& 0& 3& 1\\
    3& 5& 0& 4& 1& 2\\
    4& 3& 1& 2& 5& 0\\
    5& 1& 2& 3& 0& 4
\end{array}
\end{displaymath}
\end{example}

\begin{example}\label{Ex:E12}
An LG2-quasigroup with additional properties (see below).
\begin{displaymath}
\begin{array}{cccc}
    0& 2& 3& 1\\
    1& 3& 2& 0\\
    2& 0& 1& 3\\
    3& 1& 0& 2
\end{array}
\end{displaymath}
\end{example}

\begin{example}\label{Ex:E18}
An LG2-quasigroup with additional properties (see below).
\begin{displaymath}
\begin{array}{cccc}
    0& 3& 1& 2\\
    1& 2& 0& 3\\
    2& 1& 3& 0\\
    3& 0& 2& 1
\end{array}
\end{displaymath}
\end{example}

\begin{example}\label{Ex:E22}
An LC3-quasigroup with additional properties (see below).
\begin{displaymath}
\begin{array}{ccccc}
    0& 1& 2& 3& 4\\
    3& 2& 4& 1& 0\\
    4& 3& 1& 0& 2\\
    2& 0& 3& 4& 1\\
    1& 4& 0& 2& 3
\end{array}
\end{displaymath}
\end{example}

\begin{example}\label{Ex:E25}
An MNQ-quasigroup with additional properties (see below).
\begin{displaymath}
\begin{array}{ccccc}
    0& 1& 2& 3& 4\\
    1& 0& 3& 4& 2\\
    2& 4& 0& 1& 3\\
    3& 2& 4& 0& 1\\
    4& 3& 1& 2& 0
\end{array}
\end{displaymath}
\end{example}

\begin{example}\label{Ex:ELast}
A CQ-quasigroup with additional properties (see below).
\begin{displaymath}
\begin{array}{cccccccc}
    1&3&0&2&7&6&4&5\\
    4&2&6&0&5&3&7&1\\
    0&6&2&4&3&5&1&7\\
    7&5&4&6&1&2&0&3\\
    2&0&3&1&6&7&5&4\\
    6&4&5&7&2&1&3&0\\
    3&7&1&5&0&4&2&6\\
    5&1&7&3&4&0&6&2
\end{array}
\end{displaymath}
\end{example}

Here is the systematic search for all distinguishing examples, as outlined in
the strategy above.

$\mathcal A=$LG1: Type 1 examples are LG1$\setminus\mathcal B$, where $\mathcal
B\in\{$LC3, LNQ, LAQ, MNQ, RG2, RBQ, RG3$\}$. The remaining examples are:
$\mathcal B$=LG2 by Example \ref{Ex:E1} (check that LG2 fails with $x=0$,
$y=0$, $z=1$), $\mathcal B$=FQ by Example \ref{Ex:A3} (with $x=1$, $y=0$),
$\mathcal B$=CQ by Example \ref{Ex:A3} (with $x=0$, $y=1$, $z=0$), $\mathcal
B$=RAQ by Example \ref{Ex:A3} (with $x=0$, $y=1$), $\mathcal B$=RNQ by Example
\ref{Ex:E5} (with $x=0$, $y=0$, $z=1$), $\mathcal B$=RC3 by Example \ref{Ex:A3}
(with $x=0$, $y=0$, $z=1$), and $\mathcal B$=RC2 by Example \ref{Ex:A3} (with
$x=0$, $y=0$, $z=1$).

$\mathcal A=$LC4: Type 2 example is $\mathcal B$=LBQ. The remaining example is
$\mathcal B$=LG3 by Example \ref{Ex:E8} (with $x=0$, $y=1$, $c=0$).

$\mathcal A=$LC2: With $\mathcal B$=LC4 by Example \ref{Ex:A2} (with $x=1$,
$y=0$, $z=0$).

$\mathcal A=$LG3: With $\mathcal B\in \{$LC2, LBQ$\}$ by Example \ref{Ex:E10}
(both with $x=0$, $y=0$, $z=1$).

$\mathcal A=$LBQ: Type 2 examples are $\mathcal B\in\{$LC2, LG3$\}$.

$\mathcal A=$LG2: Type 1 examples are $\mathcal B\in\{$LC3, LNQ, LAQ, MNQ, RG2,
RBQ, RG3$\}$. The remaining examples are: $\mathcal B\in\{$LC2, LBQ, CQ$\}$ by
Example \ref{Ex:E12} (all with $x=0$, $y=0$, $z=1$), $\mathcal B$=LG3 by the
same example (with $x=0$, $y=1$, $z=0$), $\mathcal B$=FQ by Example \ref{Ex:A3}
(with $x=1$, $y=0$), $\mathcal B$=RAQ by the same example (with $x=0$, $y=1$),
$\mathcal B\in\{$RC3, RC2$\}$ by the same example (both with $x=0$, $y=0$,
$z=1$), and $\mathcal B$=RNQ by Example \ref{Ex:E18} (with $x=0$, $y=0$,
$z=1$).

$\mathcal A=$LC1: Type 2 examples are $\mathcal B\in\{$LG3, LBQ, LG2, FQ, CQ,
RAQ, RG2, RNQ, RBQ, RC3, RG3, RC2$\}$.

$\mathcal A=$LC3: Type 1 example is $\mathcal B$=MNQ. The remaining examples
are: $\mathcal B$=LC2 by Example \ref{Ex:A4} (with $x=1$, $y=0$, $z=0$),
$\mathcal B$=LAQ by the same example (with $x=1$, $y=0$), and $\mathcal B$=LNQ
by Example \ref{Ex:E22} (with $x=1$, $y=0$, $z=0$).

$\mathcal A=$LNQ: Type 1 example is $\mathcal B$=MNQ. Type 2 examples are
$\mathcal B\in\{$LC2, LC3, LAQ$\}$.

$\mathcal A=$LAQ: Type 1 example is $\mathcal B$=MNQ. Type 2 examples are
$\mathcal B\in\{$LC2, LC3, LNQ$\}$.

$\mathcal A=$EQ: Type 2 examples are $\mathcal B\in\{$LG3, LG2$\}$.

$\mathcal A=$MQ: Type 2 examples are $\mathcal B\in\{$LC2, LC3, LNQ, MNQ,
CQ$\}$.

$\mathcal A=$FQ: Type 1 examples are $\mathcal B\in\{$LBQ, LAQ$\}$.

$\mathcal A=$MNQ: Type 2 examples are $\mathcal B\in\{$LNQ, LAQ$\}$. The
remaining examples are: $\mathcal B\in\{$LC2, LC3$\}$ by Example \ref{Ex:E25}
(both with $x=1$, $y=0$, $z=2$).

$\mathcal A=$CQ: Type 1 examples are $\mathcal B\in\{$LC3, LAQ, LNQ, MNQ$\}$.
Type 2 example is $\mathcal B$=FQ. The remaining example is $\mathcal B$=LC2 by
Example \ref{Ex:ELast} (with $x=0$, $y=1$, $z=2$).

\section{Main result}

\begin{theorem} There are $26$ varieties of quasigroups of Bol-Moufang type.
Their names and defining identities can be found in Table
$\ref{Tb:Definitions}$. The maximal sets of equivalent identities for these
varieties can be found in the Hasse diagram $\ref{Fg:QTree}$, together with all
inclusions between the varieties. For every variety, we also indicate in Figure
$\ref{Fg:QTree}$ if it consists of loops, left loops or right loops.
\end{theorem}

\section{Commutative quasigroups of Bol-Moufang type}

\noindent We conclude the paper with the classification of varieties of
commutative loops of Bol-Moufang type and commutative quasigroups of
Bol-Moufang type. With the exception of the comments below, we omit all proofs
and distinguishing examples.

\begin{figure}[ht]\begin{center}\input{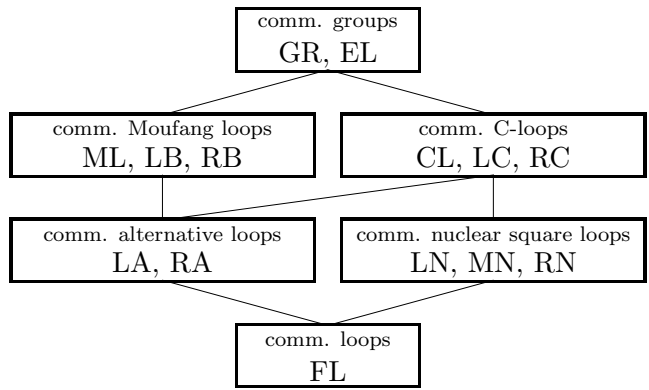}\end{center}
\caption{Varieties of commutative loops of Bol-Moufang type}\label{Fg:CL}
\end{figure}

Figure \ref{Fg:CL} depicts the Hasse diagram of commutative loops of
Bol-Moufang type. We indicate the equivalent identities defining a given
variety by pointing to the noncommutative case. For instance, commutative
groups can be defined by any defining identity of GR or EL of Figure
\ref{Fg:LTree}. Note that flexibility follows from commutativity, as $x(yx) =
x(xy) = (xy)x$.

%FIGURE
\setlength{\unitlength}{1mm}

\begin{figure}[ht]\begin{center}\input{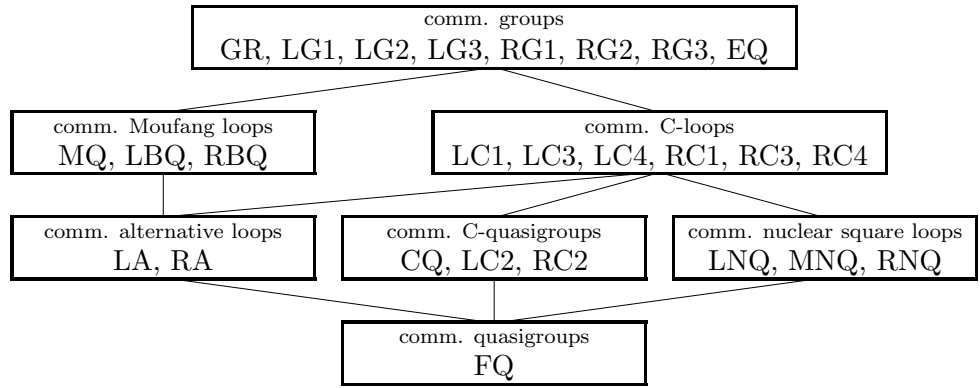}\end{center}
\caption{Varieties of commutative quasigroups of Bol-Moufang type}\label{Fg:CQ}
\end{figure}

Figure \ref{Fg:CQ} depicts the Hasse diagram of commutative quasigroups of
Bol-Moufang type, using conventions analogous to those of Figure \ref{Fg:CL}.
Note that if a quasigroup variety possesses a one-sided neutral element, its
commutative version consists of loops. Hence only the noncommutative varieties
LC2, FQ, CQ and RC2 behave differently than in the commutative loop case.
Example \ref{Ex:A2} gives a commutative C-quasigroup that is not a loop.

%FIGURE
\setlength{\unitlength}{1mm}

\section{Acknowledgement}

\noindent We thank Michael Kinyon for bringing \cite{Ku2} to our attention.

%%%%%%%%%%%%%%%%%%%%%%%%%%%%%%%%%%%%%%%%%%%%%%%%%%%%%%%%%%%%%%%%%%%%%%%%%%%%%%%
% BIBLIOGRAPHY                                                                %
%%%%%%%%%%%%%%%%%%%%%%%%%%%%%%%%%%%%%%%%%%%%%%%%%%%%%%%%%%%%%%%%%%%%%%%%%%%%%%%

\bibliographystyle{plain}

\end{document}